\input amstex
\input amsppt.sty
\input xy
\xyoption{all}


\NoBlackBoxes
\magnification=\magstephalf
\nologo

\define \Om{\Omega}
\define \n{\nabla}
\define \del{\partial}
\define \vp{\varphi}
\define \ve{\varepsilon}
\define \sdr#1#2{\overset #1\to{\underset #2\to\rightleftarrows}}
\define\tom{\widetilde \Omega}
\define \si{s\sp{-1}}

\define\sgn{\text {sgn}}
\def\Fbar{\overline F}
\def\phihat{\widehat\varphi}
\def\op{\Cal}
\def\cat{\bold}
\def\acoalg{\op A\text{-}\cat {Coalg}}
\def\fatcoalg{(\op A,\psi_{\op F})\text{-}\cat {Coalg}}
\def\aalg{\op A\text{-}\cat {Alg}}
\def\acirc{\underset \op A\to \diamond}
\def\ind{\operatorname{Ind}}
\def\tC{\widetilde C}

\topmatter
\title A canonical enriched Adams-Hilton model for simplicial sets\endtitle
\author Kathryn Hess\\ Paul-Eug\`ene Parent\\ Jonathan Scott\\ Andrew Tonks\endauthor
\date 7 September 2005\enddate
\affil Ecole Polytechnique F\'ed\'erale de Lausanne\\ University of 
Ottawa\\ Ecole Polytechnique F\'ed\'erale de Lausanne\\ London Metropolitan University \endaffil
\address Institut de g\'eom\'etrie, alg\`ebre et topologie (IGAT), Ecole Polytechnique F\'ed\'erale 
de Lausanne, CH-1015 Lausanne, Switzerland\endaddress
\email kathryn.hess\@epfl.ch\endemail
\address Department of mathematics and statistics, University of Ottawa, 585 King Edward 
Avenue, Ottawa, ON, K1N 6N5, Canada\endaddress
 \email 
pauleugene.parent\@science.uottawa.ca\endemail 
\address Institut de g\'eom\'etrie, alg\`ebre et topologie (IGAT), Ecole Polytechnique F\'ed\'erale 
de Lausanne, CH-1015 Lausanne, Switzerland\endaddress
\email jonathan.scott\@epfl.ch\endemail
\address Department of Computing, Communications Technology and Mathematics, London Metropolitan University, 166-220 Holloway Road, London N7 8DB, England\endaddress
\email a.tonks\@londonmet.ac.uk\endemail
\leftheadtext {K. Hess, P.-E. Parent, J. Scott, A. Tonks}

 \abstract For any 1-reduced simplicial set $K$ we define a canonical, 
 coassociative 
 coproduct on $\Om C(K)$, the cobar construction applied to the 
 normalized, integral chains 
 on $K$, such that any canonical quasi-isomorphism of chain algebras from 
 $\Om C(K)$ to the normalized, integral chains on $GK$, the loop group 
 of $K$, is a coalgebra map up to strong homotopy.  Our proof relies on the operadic description of the category of chain coalgebras and of strongly homotopy coalgebra maps given in \cite {HPS}. 
 \endabstract \endtopmatter
\document


\head Introduction\endhead

Let $X$ be a topological space. It is, in general, quite difficult to calculate
the algebra structure of the loop space homology 
$H\sb *\Omega X$ directly from the (singular or cubical) chain complex $C\sb *\Omega X$.  An
algorithm that associates to a space $X$ a differential graded algebra whose
homology is relatively easy to calculate and isomorphic as an algebra to
$H\sb *\Omega X$ is therefore of great value.  

In 1955 \cite {AH}, Adams and Hilton invented such an algorithm for the class of
simply-connected CW-complexes, which can be summarized as follows. 
Let $X$ be a CW-complex such that $X$ has exactly one 0-cell and no 1-cells,
and such that every attaching map is based with respect to the unique 0-cell
of $X$. There exists a morphism of differential
graded algebras inducing an isomorphism on homology---a {\sl
quasi-isomorphism}---
$$\theta \sb X: (TV,d)@>\simeq>>C\sb *\Omega X,$$
such that $\theta \sb X$ restricts to quasi-isomorphisms $(TV\sb {\leq
n},d)@>\simeq>>C\sb *\Omega X\sb {n+1}$, where $X\sb {n+1}$
denotes the $(n+1)$-skeleton of $X$, $TV$ denotes the free (tensor)
algebra on a free, graded $\Bbb Z$-module $V$,  
$\Omega X$ is the space of Moore loops on $X$ and $C_{*}$ denotes the 
cubical chains.  The morphism
$\theta\sb X$ is called an {\sl Adams-Hilton model of $X$} and satisfies the
following properties. 
\roster
\item "{$\bullet $}"If $X=\ast \cup\bigcup \sb {\alpha
\in A} e\sp {n\sb {\alpha }+1}$, then
$V$ has a degree-homogeneous basis $\{v\sb {\alpha }:\alpha \in A\}$ such
that
$\deg v\sb \alpha =n\sb {\alpha }$.
\item "{$\bullet $}"If $f\sb {\alpha }:S\sp {n\sb {\alpha }}@>>>X\sb
{n\sb{\alpha }}$ is the attaching map of the cell
$e\sp {n\sb {\alpha }+1}$, then 
$[\theta (dv\sb {\alpha })]=\Cal K\sb {n\sb {\alpha }}[f\sb {\alpha}]$.  Here,
$\Cal K\sb {n\sb {\alpha }}$ is defined so that
$$\xymatrix{
\pi \sb {n\sb {\alpha }}X\sb {n\sb {\alpha }}\ar [r]^{\cong}\ar [dr] _{\Cal K\sb {n\sb {\alpha }}}&\pi \sb
{n\sb {\alpha }-1}\Omega X\sb {n\sb {\alpha }}\ar [d]^h\\
&H\sb {n\sb {\alpha }-1}\Omega X\sb {n\sb {\alpha }}}$$
commutes, where $h$ denotes the Hurewicz homomorphism.
\endroster
It follows that $(TV,d)$ is unique up to isomorphism.

The Adams-Hilton model has proved to be a powerful tool for calculating
the loop space homology algebra of CW-complexes.  Many common spaces have
Adams-Hilton models that are relatively simple and thus well-adapted to
computations.  Difficulties start to arise, however, when one wishes to use the Adams-Hilton model 
to compute the algebra homomorphism induced by a cellular map 
$f:X@>>>Y$.  

If $\theta \sb X:(TV,d)@>>>C\sb *\Omega X$ and
$\theta \sb Y:(TW,d)@>>>C\sb *\Omega Y$ are Adams-Hilton models, then
there exists a unique homotopy class of morphisms $\vp :
(TV,d)@>>>(TW,d)$ such that 
$$\xymatrix{
(TV,d)\ar [r]^{\vp}\ar [d] ^{\theta _{X}}&(TW,d)\ar [d] ^{\theta 
_{Y}}\\
C\sb *\Omega X\ar [r] ^{C_{*}\Om f}&C\sb *\Omega Y}$$
commutes up to derivation homotopy.¥
Any representative $\vp$ of this
homotopy class can be said to be an Adams-Hilton model of $f$.   

As the choice of $\vp $ is unique only up to homotopy, the Adams-Hilton 
model is not a functor.  The essential problem is that choices are 
made at each stage of the construction of  
$\theta _{X}$ and $\theta _{Y}$: they are not canonical. For many purposes this lack of 
functoriality does not cause any problems.  When one needs to use 
Adams-Hilton models to construct new models, however, then it can become quite troublesome, as seen in, 
e.g., \cite {DH} .  

Similarly, when constructing algebraic models based on Adams-Hilton 
models, one often needs the models to be {\sl enriched}, i.e., there should 
be a chain algebra map $\psi: (TV,d)@>>>(TV,d)\otimes (TV,d)$ such 
that  
$$\xymatrix{
(TV,d)\ar [rrr]^{\psi}\ar [d] ^{\theta _{X}}&&&(TV,d)\otimes (TV,d)\ar [d] ^{\theta 
_{X}\otimes \theta _{X}¥}\\
C\sb *\Omega X\ar [rrr] ^{AW\circ C_{*}\Om \Delta _{X}¥}&&&C\sb *\Omega X\otimes C\sb *\Omega X }$$
commutes up to homotopy, where $AW$ denotes the Alexander-Whitney 
equivalence.  Thus $\theta _{X}$ is a coalgebra map up to homotopy.  Since the underlying algebra of
$(TV,d)$ is free, such a coproduct 
always exists and can be constructed degree by degree.  Again, however, choices are involved in the 
construction of $\psi $, so that one usually knows little about it, other 
than that it exists. In particular, since the diagram above 
commutes and $AW\circ C_{*}\Om \Delta _{X}$ is cocommutative up to 
homotopy and strictly coassociative, $\psi$ 
is coassociative and cocommutative up to homotopy, i.e., 
$(TV,d,\psi)$ is a {\sl Hopf algebra up to homotopy} \cite {A}.   For many 
constructions, however, it would be very helpful to know that there 
is a choice of $\psi$ that is strictly coassociative.

Motivated by the need to rigidify the Adams-Hilton model 
construction and its enrichment, we work here with simplicial sets rather than 
topological spaces.  Any topological space $X$ that is equivalent to a finite-type 
simplicial complex is homotopy-equivalent to the realization of a 
finite-type simplicial set.  There is an obvious candidate for a canonical 
Adams-Hilton model of a $1$-reduced simplicial set $K$:    $\Om C(K)$, the 
cobar construction on the integral, normalized chains on $K$, which is 
a free algebra on generators in one-to-one correspondence with the nondegenerate 
simplices of $K$. It follows easily by acyclic models methods (see, 
e.g, \cite {M}) that there exists a natural quasi-isomorphism of chain 
algebras $$\theta_{K}: \Om C(K)@>\simeq>> C(GK),$$ 
where $GK$ denotes the Kan loop group on $K$. There is also an explicit formula for such a natural transformation, due to Szczarba \cite 
{Sz}.

In this article we provide a simple definition of a natural, strictly 
coassociative coproduct, the {\sl Alexander-Whitney (A-W) cobar 
diagonal},
$$\psi _{K}: \Om C(K)@>>>\Om C(K)\otimes \Om C(K),$$
where $K$ is any 1-reduced simplicial set.
Furthermore, any natural quasi-isomorphism of chain algebras
$\theta_{K}: \Om C(K)@>\simeq>> C(GK)$ is a strongly homotopy coalgebra map with respect to $\psi _{K}$.  In other words, $\theta_{K}$ 
is a coalgebra map up to homotopy; the homotopy in question is a 
coderivation up to a second homotopy; etc. The map $\theta_{K}$ has already 
proved extremely useful in 
constructing a number of interesting algebraic models, such as in \cite 
{BH},\cite {H}, \cite 
{HL}. 

Ours is not the only definition of a canonical, coassociative coproduct on $\Om 
C(K)$. In \cite {Ba} Baues defined combinatorially an explicit coassociative 
coproduct $\widetilde \psi _{K}$¥ on $\Om 
C(K)$, together with an explicit derivation homotopy insuring 
cocommutativity up to homotopy.  He showed that there is an injective 
quasi-isomorphism of chain Hopf algebras from $(\Om C(K), \widetilde 
\psi _{K})$ into (the first Eilenberg subcomplex of) the cubical 
cochains on the geometrical cobar construction on $K$.  

We show in section 5 of this article that Baues's coproduct is equal to the Alexander-Whitney cobar diagonal, a result that is surprising at first sight. It is clear from the definition of Baues's coproduct that its image lies in $\Om C(K)\otimes \si C_+(K)$, so that its form is highly asymmetric.  That asymmetry is well hidden in our definition of the Alexander-Whitney cobar diagonal.

Even though the two definitions are equivalent, our 
approach is still interesting, as the Alexander-Whitney cobar diagonal is 
given explicitly in terms of only two fundamental pieces: the 
diagonal map on a simplicial set and the Eilenberg-Zilber strong 
deformation retract
$$C(K)\otimes C(K)\sdr {\n }{f}C(K\times K)\circlearrowleft\vp.
$$
(See section 2.)  Furthermore, it is very helpful for construction 
purposes to have an explicit equivalence $\theta_{K}:\Om 
C(K)@>>>C(GK)$ that is a map of coalgebras up to strong homotopy and a map of algebras, as the articles \cite 
{BH},\cite {H}, \cite 
{HL} amply illustrate. 

In a subsequent article \cite {HPS2}, we will further demonstrate the importance of our coproduct definition on the cobar construction, when we treat the special case of suspensions.  In particular we will show that the Szczarba equivalence  is a strict coalgebra map when $K$ is a suspension.

After recalling a number of elementary definitions at the end of this 
introduction, we devote section 1 to Gugenheim and Munkholm's category $\bold{DCSH}$, the category of chain coalgebras and strongly homotopy coalgebra maps.  In particular, we recall and expand upon the ``operadic" description of $\bold{DCSH}$ developed in \cite{HPS}. In section 2 we introduce homological 
perturbation theory and its interaction with morphisms in $\bold{DCSH}$, 
expanding the discussion to include twisting cochains 
and twisting functions in section 3.  The heart of this article is section 4, where we define the Alexander-Whitney cobar 
diagonal, show that it is cocommutative up to homotopy and strictly 
coassociative, and prove that the Szczarba equivalence  
is a strongly homotopy coalgebra map. We conclude section 4 with a discussion of the relationship of our work to the problem of iterating the cobar construction. In section 5 we prove that the Alexander-Whitney cobar diagonal is equal to Baues's coproduct on $\Om C(K)$. 

In a forthcoming paper we will explain how the canonical 
Adams-Hilton model enables us to carry out Bockstein spectral sequence calculations using methods previously applied only in the ``Anick" range (cf., \cite {Sc}) to spaces well outside of that range.

\subhead Preliminary definitions, terminology and notation\endsubhead

We recall here certain necessary elementary definitions and constructions.  We also introduce notation 
and terminology that we use throughout the remainder of this paper.

We consider that the set of natural numbers $\Bbb N$ includes 0.

If $\bold C$ is a category and $A$ and $B$ are objects in $\bold C$, then $\bold C(A,B)$ denotes the collection of morphisms from $A$ to $B$. We write $\bold C^\rightarrow$ for the category of morphisms in $\bold C$.

Given chain complexes $(V,d)$ and $(W,d)$, the notation
$f:(V,d)@>\simeq >>(W,d)$ indicates that $f$ induces an isomorphism in homology. 
In this case we refer to $f$ as a {\sl quasi-isomorphism}.

The {\sl suspension} endofunctor $s$ on the category of graded modules is defined on objects $V=\bigoplus \sb {i\in \Bbb Z} V\sb i$ by
$(sV)\sb i \cong V\sb {i-1}$.  Given a homogeneous element $v$ in
$V$, we write $sv$ for the corresponding element of $sV$. The suspension $s$ admits an obvious inverse, which we denote $\si$.

A graded $R$-module $V=\bigoplus \sb {i\in \Bbb Z} V\sb i$ is {\sl connected} 
if $V_{<0}=0$ and $V_{0}\cong R$.  It is {\sl simply connected} if, in 
addition, $V_{1}=0$.  We write $V_{+}$ for $V_{>0}$.

Let $V$ be a positively-graded $R$-module.  The free associative algebra on $V$ is denoted 
$TV$, i.e., 
$$TV\cong R\oplus V\oplus (V\otimes V)\oplus (V\otimes V\otimes V)\oplus\cdots .$$
A typical basis element of $TV$ is denoted $v\sb 1\cdots v\sb n$, i.e., we drop the
tensors from the notation.  We say that $v_1\cdots v_n$ is {\sl of length $n$} and let $T^nV=V^{\otimes n}$ be the submodule of words of length $n$. The product on $TV$ is then defined by
$$\mu (u\sb 1\cdots u\sb m\otimes v\sb 1\cdots v\sb n)=u\sb 1\cdots u\sb m v\sb
1\cdots v\sb n.$$
Throughout this paper $\pi:T^{>0}V@>>>V$ denotes the projection map such that  $\pi (v_1\cdots v_n)=0$ if $n>1$ and $\pi (v)=v$ for all $v\in V$.  When we refer to the {\sl linear part} of an algebra map $f:TV@>>>TW$, we mean the composite $\pi \circ f|_V:V@>>>W$.

\definition {Definition} Let $(C,d)$ be a simply-connected chain 
coalgebra with reduced coproduct $\overline \Delta$.  The {\sl cobar construction} on $(C,d)$, denoted $\Om 
(C,d)$, is the chain algebra $(T\si (C_{+}),d_{\Om})$, where 
$d_{\Om}=-\si ds+(\si \otimes \si)\overline \Delta s$ on generators.\enddefinition

Observe that for every pair of simply-connected chain coalgebras $(C,d)$ and $(C',d')$ 
there is a natural quasi-isomorphism of chain algebras
$$q:\Om \bigl( (C,d)\otimes (C',d')\bigr )@>>>\Om (C,d)\otimes \Om 
(C',d')\tag 0.1$$
specified by $q\bigl( \si (x\otimes 1)\bigr)=\si x\otimes 1$, $q\bigl( \si (1\otimes 
y)\bigr)=1\otimes\si y$ and $q\bigl( \si (x\otimes y)\bigr)=0$ for all $x\in 
C_{+}$ and $y\in C'_{+}$ \cite {Mi:Thm. 7.4}.

\definition {Definition}  Let $f,g:(A,d)@>>>(B,d)$ be two maps of 
chain algebras.  An {\sl $(f,g)$-derivation} is a linear map $\vp 
:A@>>>B$ of degree $+1$ such that  $\vp \mu 
=\mu (\vp \otimes g+f \otimes \vp)$, where $\mu$ denotes the 
multiplication on $A$ and $B$.   A {\sl derivation homotopy} from $f$ to $g$ is an $(f,g)$-derivation $\vp $that satisfies $d\varphi + \varphi d = f - g$. \enddefinition

If $f$ and $g$ are maps of chain coalgebras, there is an obvious 
dual definition of an {\sl $(f,g)$-coderivation} and of {\sl $(f,g)$-coderivation homotopy}.

\definition {Definition} Let $K$ be a simplicial set, and let $\Cal 
F_{ab}$ denote the free abelian group functor.  For all $n>0$, let 
$DK_{n}=\cup _{i=0}^{n-1} s_{i}(K_{n-1})$, the set of degenerate 
$n$-simplices of $K$. The {\sl 
normalized chain complex} on $K$, denoted $C(K)$, is given by
$$C_{n}(K)=\Cal F_{ab} (K_{n})/\Cal F_{ab}(DK_{n}).$$
Given a map of simplicial sets $f:K@>>>L$, the induced map of 
normalized chain complexes is denoted $f_{\sharp}$.
\enddefinition

\definition {Definition} Let $K$ be a reduced simplicial set, and let $\Cal F$ denote the 
free group functor. The {\sl loop group} $GK$ on $K$ is the simplicial group such that $(GK)_{n}=\Cal F 
(K_{n+1}\smallsetminus \operatorname {Im} s_{0})$, with faces and degeneracies 
specified by
$$\split 
\del _{0}\bar x&=(\;\overline {\del _{0}x}\;)^{-1}\overline {\del 
_{1}x}\\
\del _{i}\bar x&=\overline {\del 
_{i+1}x}\quad\text {for all $i>0$}\\
s_{i}\bar x&=\overline {s_{i+1}x}\quad \text {for all $i\geq 0$}
\endsplit$$
where $\bar x$ denotes the class in $(GK)_{n}$ of $x\in K_{n+1}$.
\enddefinition

Observe that for each pair of reduced simplicial sets $(K,L)$ there is a 
unique homomorphism of simplicial groups $\rho :G(K\times L)@>>>GK\times 
GL$, which is specified by $\rho \bigl(\;\overline {(x,y)}\;\bigr)=(\bar x, \bar y)$.


\head 1. The category $\bold {DCSH}$ and its relatives\endhead

The categor $\bold {DCSH}$ of coassociative chain coalgebras and of coalgebra morphisms up to strong homotopy was first defined by Gugenheim and Munkholm in the early 1970's \cite {GM}, when they were studying extended naturality of the functor $\operatorname{Cotor}$. The objects of $\bold {DCSH}$  have a relatively simple algebraic description, while that of the morphisms is rich and complex.  Its objects are augmented, coassociative chain coalgebras, and a morphism from $C$ to $C'$ is a map of chain algebras $\Omega C@>>>\Omega C'$.  

In a slight abuse of terminology, we say that a chain map between chain coalgebras $f:C@>>>C'$ is a {\sl DCSH map} if there is a morphism  in $\bold {DCSH}(C,C')$ of which $f$ is the linear part. In other words, there is a map of chain algebras $g:\Om C@>>>\Om C'$ such that 
$$g|_{s^{-1}C_+}=\si f s +\text{higher-order terms}.$$
In a further abuse of notation, we sometimes write $\tom f:\Om C@>>>\Om C'$ to indicate one choice of chain algebra map of which $f$ is the linear part.

It is also possible to broaden the definition of coderivation 
homotopy to homotopy of DCSH maps.  Given two DCSH maps $f,f ':C@>>>C'$, a {\sl DCSH homotopy} from $f$ to $f'$  is a  $(\tom f,\tom f ')$-derivation 
homotopy $h:\Om (C,d)@>>>\Om (C',d')$.  We sometimes abuse terminology and refer to the linear part of $h$ as a DCSH homotopy from $f$ to $f'$.    

The category $\bold {DCSH}$ plays an important role in topology.  For any reduced simplicial set $K$, the usual coproduct on $C(K)$  is a DCSH map, as we explain in detail in section 2.  Furthermore, we show in section 4 that given any natural, strictly coassociative coproduct on $\Om C(K)$, any natural map of chain algebras $\Om C(K)@>>>C(GK)$ is also a DCSH map.

In \cite {HPS} the authors provided a purely operadic description of $\bold {DCSH}$.  Before recalling and elaborating upon this description, we briefly explain the framework in which it is constructed.  We refer the reader to section 2 of \cite {HPS} for further details. 

Let $\bold M$ denote the category of chain complexes over a PID $R$, and let $\bold M^\Sigma $denote the category of symmetric sequences of chain complexes.  An object $\Cal X$ of $\bold M^\Sigma$ is  a family $\{\Cal X(n)\in \bold M\mid n\geq 0\}$ of objects in $\bold M$ such that $\Cal X(n)$ admits a right action of the symmetric group $\Sigma_n$, for all $n$.  There is a faithful functor $\xymatrix@1{\Cal T: \bold M\ar [r]&\bold M^\Sigma}$ where, for all $n$, $\Cal T(A)(n)=A^{\otimes n}$, where $\Sigma _n$ acts by permuting the tensor factors.  The functor $\Cal T$ is strong monoidal, with respect to the {\sl level monoidal structure} $(\bold M^\Sigma, \otimes, \op C)$, where $(\op X\otimes \op Y)(n)=\op X(n)\otimes\op Y(n)$, endowed with the diagonal action of $\Sigma_n$, and $\op C(n)=R$, endowed with the trivial $\Sigma _n$-action.

The category $\bold M^\Sigma$ also admits a nonsymmetric, right-closed monoidal structure $(\bold M^\Sigma,\diamond , \Cal J)$, where $\diamond$ is the {\sl composition product} of symmetric sequences, and $\Cal J(1)=R$ and $\Cal J(n)=0$ otherwise.  Given symmetric sequences $\op X$ and $\op Y$, $(\op X\diamond \op Y)(0) =\op X(0)\otimes \op Y(0)$ and for $n>0$,
$$(\op X\diamond\op Y)(n)=\coprod\Sb k\geq 1\\ \vec\imath\in I_{k,n}\endSb \op X(k)\underset \Sigma _k\to \otimes \bigl(Y(i_1)\otimes \cdots\otimes Y(i_k)\bigr)\underset \Sigma _{\vec\imath}\to \otimes R[\Sigma _n],$$
where $I_{k,n}=\{\vec\imath=(i_1,...,i_k)\in \Bbb N^k\mid \sum _j i_j=n\}$ and $\Sigma_{\vec\imath}=\Sigma _{i_1}\times \cdots\times \Sigma _{i_k}$, seen as a subgroup of $\Sigma _n$.
For any objects $\Cal X, \Cal X', \Cal Y, \Cal Y'$ in $\bold M^\Sigma$, there is an obvious, natural intertwining map
$$\xymatrix{\iota: (\Cal X\otimes \Cal X')\diamond (\Cal Y\otimes \Cal Y')\ar [r]&(\Cal X\diamond \Cal Y)\otimes (\Cal X'\diamond\Cal Y')}.\tag 2.1$$

An {\sl operad} in $\bold M$ is a monoid with respect to the composition product.  The {\sl associative operad} $\Cal A$ is given by $\Cal A(n)=R[\Sigma _n]$ for all $n$, endowed with the obvious monoidal structure, induced by permutation of blocks.

Let $\Cal P$ denote any operad in $\bold M$.  A {\sl $\Cal P$-coalgebra} consists of an object $C$ in $\bold M$. together with an appropriately equivariant and associative family $$\{\xymatrix@1{C\otimes \Cal P (n) \ar [r]&C^{\otimes n}\mid n\geq 0}\}$$ of morphisms in $\bold M$.  The functor $\Cal T$ restricts to a faithful functor $$\xymatrix@1{\Cal T: \op P\text{-}\cat {Coalg}\ar [r]&\bold{Mod}_{\Cal P}}$$ from the category of $\Cal P$-coalgebras to the category of right $\Cal P$-modules.

In \cite {HPS} the authors constructed a free $\Cal A$-bimodule $\Cal F$, called the {\sl Alexander-Whitney bimodule}.  As symmetric sequences of graded modules, $\Cal F =\op A\diamond \op S\diamond \op A$, where $\op S(n)= R[\Sigma _n]\cdot z_{n-1}$, the free $R[\Sigma _n]$-module on a generator of degree $n-1$. Moreover, $\op F$ admits an increasing, differential filtration, given by $F_n \op F=\op A \diamond\op S_n\diamond \op A$, where $\op S_n(m)=\op S (m)$ if $m\leq n$ and $\op S_n(m)=0$ otherwise.  More precisely, if $\del _{\op F}$ is the differential on $\op F$, then
$$ \del _{\op F} z_n=\sum _{0\leq i\leq n-1} \delta \otimes (z_i\otimes z_{n-i-1}) +\sum _{0\leq i\leq n-1} z_{n-1}\otimes (1^{\otimes i}\otimes \delta\otimes 1^{\otimes n-i-1}),$$
where $\delta\in \op A(2)=R[\Sigma _2]$ is a generator.

The Alexander-Whitney bimodule is endowed with a coassociative, counital coproduct 
$$\xymatrix@1 {\psi_{\op F}: \Cal F\ar [r]&\Cal F\acirc \Cal F},$$
where $ \acirc $ denotes the composition product over $\Cal A$, defined as the obvious coequalizer.  In particular,
$$\psi_{\op F} (z_n)=\sum\Sb 1\leq k\leq n+1\\ \vec \imath\in I_{k,n+1}\endSb z_{k-1}\otimes (z_{i_1-1}\otimes \cdots \otimes z_{i_k-1})$$
for all $n\geq 0$, where $I_{k,n}=\{\vec\imath=(i_1,...,i_k)\mid \sum _j i_j=n\}$.   

Furthermore, $\Cal F$ is a level comonoid, i.e., there is a coassociative, counital coproduct $$\xymatrix@1{\Delta _{\op F}:\op F\ar [r]&\op F\otimes \op F},$$
which is specified by 
$$\Delta_{\op F}(z_n)=\sum\Sb 1\leq k\leq n+1\\ \vec \imath\in I_{k,n+1}\endSb\bigl(z_{k-1}\otimes (\delta ^{(i_1)}\otimes \cdots \otimes \delta ^{(i_k)})\bigr)\otimes \bigl (\delta ^{(k)}\otimes (z_{i_1 -1}\otimes \cdots \otimes z_{i_k-1})\bigr).$$
Here, $\delta ^{(i)}\in \op A(i)$ denotes the appropriate iterated composition product of $\delta ^{(2)}=\delta$.

Let $(\Cal A,\psi_{\op F})\text{-}\bold {Coalg}$ denote the category of which the objects are $\Cal A$-coalgebras (i.e., coassociative and counital chain coalgebras)  and where the morphisms are defined by 
$$\fatcoalg(C,C'):=\bold{Mod}_{\Cal A}\bigl(\Cal T(C)\underset \Cal A\to \circ\Cal F, \Cal T(C')\bigr).$$
Composition in $(\Cal A,\psi_{\op F})\text{-}\bold {Coalg}$ is defined in terms of $\psi_{\op F}$. Given $\theta\in (\Cal A,\psi_{\op F})\text{-}\bold {Coalg}(C, C')$ and $\theta'\in (\Cal A,\psi_{\op F})\text{-}\bold {Coalg}(C', C'')$, their composite $\theta'\theta\in(\Cal A,\psi_{\op F})\text{-}\bold {Coalg}(C, C'')$ is given by composing the following sequence of (strict) morphisms of right $\op A$-modules.
$$\xymatrix@1{\op T(C)\acirc \op F\ar [rr]^{1_{\op T(C)}\acirc \psi_{\op F}}&&\op T(C)\acirc \op F \acirc \op F\ar [rr]^{\theta \acirc 1_{\op F}}&&\op T(C')\acirc \op F\ar [r]^{\theta'}&\op T(C'')}.$$
We call $(\Cal A,\psi_{\op F})\text{-}\bold {Coalg}$ the {\sl $(\op F, \psi_{\op F})$-governed category of $\op A$-coalgebras}.

The  important properties of the Alexander-Whitney bimodule given below follow immediately from the Cobar Duality Theorem in \cite {HPS}. 

\proclaim{Theorem 1.1\cite{HPS}} 
For any category $\bold D$, there is a full and faithful functor, called the {\sl induction functor},
$$\ind : \bigl (\fatcoalg\bigr)^{\bold D}@>>>\bigl (\aalg\bigr)^{\bold D}$$
defined on objects by $\ind (X)=\Om X$ for all functors $X:\bold D@>>>\fatcoalg$ and on morphisms by
$$\ind(\tau)|_{\si X}=\sum _{k\geq 1} (\si )^{\otimes k} \tau (-\otimes z_{k-1} ) s :\si X @>>> \Om Y$$
for all natural transformations $\tau :X@>>> Y$.
\endproclaim

As an easy consequence of Theorem 1.1, we obtain the following result.

\proclaim {Corollary 1.2 \cite {HPS}} The category  $\cat{DCSH}$, is isomorphic to the $(\op F, \psi_{\op F})$-governed category of coalgebras, $(\Cal A,\psi_{\op F})\text{-}\cat {Coalg}$.\endproclaim

The isomorphism of the corollary above is given by the identity on objects and $\ind$ on morphisms.

Define a bifunctor $\wedge :\fatcoalg\times \fatcoalg @>>> \fatcoalg$ on objects by $C\wedge C':= C\otimes C'$, the usual tensor product of chain coalgebras.  Given $ \theta\in \fatcoalg (C, D)$ and $\theta'\in \fatcoalg (C',D')$, we define $\theta \wedge \theta'$ to be the composite of (strict) right $\op A$-module maps
$$\xymatrix{\op T(C\wedge C')\acirc \op F\ar [r]^(0.4)\cong\ar [rrddd]_{\theta \wedge \theta '}&\bigl(\op T(C)\otimes \op T(C')\bigr)\acirc \op F\ar [r]^(0.45){1\acirc \Delta _{\op F}}&\bigl(\op T(C)\otimes \op T(C')\bigr)\acirc (\op F\otimes \op F)\ar [d]^{\iota}\\
&&\bigl(\op T(C)\acirc \op F\bigr )\otimes\bigl(\op T(C')\acirc \op F\bigr )\ar [d]^{\theta \otimes \theta '}\\
&&\op T(D)\otimes \op T(D')\ar [d]^\cong\\
&&\op T(D\wedge D')}$$
where $\iota $ is the intertwining map of (2.1).
It is straightforward to show that $\wedge$ endows $\fatcoalg$ with the structure of a monoidal category.

\proclaim{Lemma 1.3}  The induction functor $\ind : \fatcoalg@>>>\aalg$ is comonoidal.\endproclaim

\demo{Proof}  Let $q:\Om (-\otimes -)@>>> \Om (-)\otimes \Om (-)$ denote Milgram's natural transformation (0.1) of functors from $\acoalg$ into $\aalg$.  It is an easy exercise, based on the explicit formula for $\Delta _{\op F}$, to prove that 
$$q\ind (\theta\wedge\theta ')=\bigl(\ind(\theta)\otimes \ind(\theta')\bigr)q:\Om (C\otimes C')@>>>\Om D\otimes \Om D'$$
for all  $ \theta\in \fatcoalg (C, D)$ and $\theta'\in \fatcoalg (C',D')$.
Milgram's equivalence therefore provides us with the desired natural transformation
$$q: \ind (-\wedge -)@>>>\ind (-)\otimes \ind (-).\qed$$
\enddemo

In section 4 of this article we consider objects in 
the following category related to $\fatcoalg$.

\definition{Definition}  The objects of the {\sl weak Alexander-Whitney category} $\bold {wF}$ are pairs $(C,\Psi)$, where $C$ is a object in $\acoalg$ and $\Psi\in \fatcoalg(C, C\otimes C)$ such that
$$\Psi (-\otimes z_0):C@>>>C\otimes C$$
is exactly the coproduct on $C$, while
$$\bold F\bigl( (C,\Psi), (C',\Psi')\bigr)= \{ \theta \in \fatcoalg(C,C')\mid \Psi' \theta =(\theta\wedge \theta )\Psi\}.$$ An object of $\cat {wF}$ is called a {\sl weak Alexander-Whitney coalgebra}. 
\enddefinition

As we establish in the next lemma, the cobar construction provides an important link between the weak Alexander-Whitney category and the following category of algebras endowed with coproducts.

\definition{Definition} The objects of the {\sl weak Hopf algebra category} $\bold {wH}$ are pairs $(A,\psi)$, where $A$ is a chain algebra over $R$ and $\psi:A@>>>A\otimes A$ is a map of chain algebras, while
$$\bold{wH}\bigl((A,\psi), (A',\psi ')\bigr)=\{ f\in \aalg(A,A')\mid \psi' f=(f\otimes f)\psi\}.$$
\enddefinition

\proclaim {Lemma 1.4} The cobar construction extends to a functor $\tom : \bold {wF}@>>>\bold {wH}$.
\endproclaim 

\demo{Proof}  Given an object $(C,\Psi)$ of $\bold {wF}$, let $\tom (C,\Psi)=\bigl(\Om C, q \ind(\Psi)\bigr)$, where $\ind(\Psi):\Om C@>>>\Om(C\otimes C)$, as in Theorem 1.1(i), and $q:\Om (C\otimes C)@>>>\Om C\otimes \Om C$ is Milgram's equivalence (0.1).  In particular, $q\ind(\Psi):\Om C@>>>\Om C\otimes \Om C$ is indeed a morphism of algebras, as it is a composite of two algebra maps.

On the other hand, given $\theta \in \bold {wF}\bigl( (C,\Psi), (C',\Psi')\bigr)$, let 
$\tom \theta =\ind (\theta):\Om C@>>>\Om C'$.  Then
$$\align
\bigl(q\ind(\Psi')\bigr)\tom\theta=&q\ind(\Psi')\ind(\theta)=q\ind(\Psi'\theta)\\
=&q\ind \bigl((\theta \wedge\theta)\Psi\bigr)=q\ind(\theta\wedge\theta)\ind(\Psi)\\
=&\bigl(\ind(\theta)\otimes \ind(\theta\bigr)\bigr )q\ind(\Psi)\\
=&\bigl(\tom(\theta)\otimes\tom(\theta)\bigr)\bigl(q\ind(\Psi)\bigr),\endalign$$
i.e., $\tom\theta$ is indeed a morphism in $\bold{wH}$.
\qed\enddemo

We are, of course, particularly interested in those objects $(C, \Psi)$ of $\cat {wF}$ for which $\tom (C, \Psi)$ is actually a strict Hopf algebra, i.e., such that $q\ind(\Psi)$ is coassociative.

\definition{Definition} The {\sl Alexander-Whitney category} $\cat F$ is the full subcategory of $\cat {wF}$ such that $(C,\Psi)$ is an object of $\cat F$ if and only if $q\ind(\Psi)$ is coassociative.  The objects of $\cat F$ are called {\sl Alexander-Whitney coalgebras}.\enddefinition

As we explain in section 4, for any reduced simplicial set $K$, there is a canonical choice of $\Psi_K$ such that $(C(K), \Psi _K)$ is an object of $\bold F$.

From the proof of Lemma 1.4, it is clear that $\tom$ restricts to a functor $\tom:\cat F@>>>\cat H$, where $\cat H$ is the category of Hopf algebras. 

\proclaim{Theorem 1.5}  Let $X,Y: \cat D@>>>\cat {H}$ be functors, where $\cat D$ is a category admitting a set of models $\frak M$ with respect to which  $Y$ is acyclic. Suppose that $X$ factors through $\bold F$ as follows
$$\xymatrix{\cat D\ar [rr]^X\ar [dr]_C&&\cat {H}\\ &\bold F\ar [ur]_\tom}$$
where $C$ is free with respect to $\frak M$. Let $\theta : UX@>>>UY$ be any natural transformation of functors into $\aalg$, where $U:\bold {H}@>>>\aalg$ denotes the forgetful functor.  Then there exists a natural transformation $\widehat \theta :\Om X@>>>\Om Y$ extending the desuspension of $\theta$, i.e., for all objects $d$ in $D$,
$$\widehat\theta (d)=\si \theta (d) s + \text{higher-order terms}.$$
\endproclaim

The proof of this result depends strongly on the notion of a free functor with respect to a set of models, which we recall  in detail, before commencing the proof of the theorem.  

Let $\cat D$ be a  category, and let $\frak M$ be a set of objects in $\cat D$.  A functor $X: \cat D@>>>\cat M$ is {\sl free} with respect to $\frak M$ if there is a set $\{e_M\in X(M)\mid M\in \frak M\}$ such that $\{X(f)(e_M)\mid f\in \cat D(M,D), M\in \frak M\}$ is an $R$-basis of $X(D)$ for all objects $D$ in $\cat D$.  If $X:\cat D@>>>\bold {H}$, then $X$ is {\sl free}Ê with respect to $\frak M$ if $U'X:\cat D@>>>\cat M$ is free, where $U':\bold {H}@>>>\cat M$ is the forgetful functor.  

\demo{Proof} According to Theorem 1.1, it suffices to construct a natural transformation
$\tau : \op T(C)\underset \op A\to \circ \op F@>>>\op T(Y)$ of right $\op A$-modules such that $\tau (-\otimes z_0)=\theta$.  We can then set $\widehat \theta =\operatorname {Ind}(\tau)$. 

Since $ \op T(C)\underset \op A\to \circ \op F=\op T(C)\circ \op S\circ \op A$, which is a free right $\op A$-module, any natural tranformation of functors into the category of symmetric sequences $\op T(C)\circ \op S@>>>\op T(Y)$ can be freely extended to a natural tranformation of functors into the category of right $\op A$-modules.  Furthermore, any family of equivariant natural transformations of functors into the category of {\sl graded $R$-modules}  
$$\{\tau _k: C\otimes \op S(k)@>>>Y^{\otimes k}\mid k\geq 1\}\tag 1.2$$
 induces a natural tranformation of functors into symmetric sequences of {\sl graded modules} $\tau': \op T(C)\circ \op S@>>>\op T(Y)$, given by composites like
$$\xymatrix {\op T(C)(k)\otimes \op S(n_1)\otimes \cdots \otimes \op S(n_k)\ar[r]^\cong& \bigl(C\otimes \op S(n_1)\bigr)\otimes \cdots \otimes \bigl(C\otimes \op S(n_k)\bigr)\ar[d]^{\tau _{n_1}\otimes\cdots\otimes \tau _{n_k}}\\
&Y^{\otimes n_1}\otimes\cdots \otimes Y^{\otimes n_k}\ar [d]^=\\
&Y^{\otimes n},}$$
where $n=\sum _i n_i$.  The free extension of $\tau '$ to  $\tau : \op T(C)\underset \op A\to \circ \op F@>>>\op T(Y)$ will be a natural transformation of functors into $\bold M$, the category of chain complexes, if  for all $k$,
$$\del_{Y^{\otimes k}}\tau _k=\tau '(\del _C\otimes 1) +\tau (1\otimes \del_{\op F})\tag 1.3$$
where $\del _{Y^{\otimes k}}$, $\del _C$ and $\del _{\op F}$ are the differentials on $Y^{\otimes k}$, $C$, and $\op F$, respectively.   Note that the formula for the restriction of $\tau'$ to $\bigl(\op T(C)\circ \op S\bigr)(k)$ involves only the $\tau _j$'s for $j\leq k$, as does the restriction of $\tau$ to $\bigl(\op T(C)\circ \op S\circ \op A\bigr)(k)$.

We construct the family (1.2) recursively.  We can choose $\tau _1$ to be the ``linear part" of $\theta$, i.e., for all objects $d$ in $\cat D$, the map of graded modules $\tau_1 (d)(-\otimes z_0)$ is the composite
$$C(d)@>\si >>\si C(d)@>\theta (d)|_{\si C(d)}>> \Om Y(d) @>\pi >> \si Y(d)@>s>> Y(d).$$
Assume that $\tau _k$ has been defined  for all $k<n$ and that $\tau _n(e_M\otimes z_{n-1})$ has been defined for all $M$ such that $\deg e_M<m$ so that (1.3) holds.  Let $M$ be an element of $\frak M$ such that $\deg e_M=m$.    According to the induction hypothesis, both $\tau ' (\del _{C(M)} e_M\otimes z_{n-1})$ and $\tau (e_M\otimes \del _{\op F} z_{n-1})$ have already been defined. Moreover, since $\del _{C(M)} e_M\otimes z_{n-1}+ (-1)^n e_M\otimes \del _{\op F} z_{n-1}$ is a cycle,   $$\tau '(\del _{C(M)} e_M\otimes z_{n-1})+(-1)^n\tau (e_M\otimes \del _{\op F} z_{n-1})$$ is a cycle in $Y(M)$ and therefore a boundary, since $Y$ is acyclic with respect to $M$.  We can thus continue the recursive construction of $\tau _n$.
\qed\enddemo


\head 2. Homological perturbation theory\endhead

In this section we recall those elements of homological perturbation 
theory that we use in the construction of the Alexander-Whitney cobar 
diagonal.

 \definition {Definition} Suppose that $\n :(X, \del ) @>>> (Y,d)$ and
$f:(Y,d)@>>>(X,\del )$ are morphisms of chain
complexes.  If $f\n = 1_{X}¥$  and there exists a
chain homotopy
$\vp : (Y,d)@>>> (Y,d)$ such that 
\roster
\item $d\vp +\vp d =\n f -1_{Y}¥$,
\item $\vp \n =0$,
\item $f\vp =0$, and
\item $\vp \sp 2=0$,
\endroster
then $(X,d) \sdr{\n}f (Y,d)\circlearrowleft\vp$ is a {\sl strong
deformation retract (SDR) of chain complexes.}  
\enddefinition

It is easy to show that given a chain homotopy $\vp '$
satisfying condition (1), there exists a chain homotopy $\vp $ satisfying
all four conditions.  As explained in, e.g.,
\cite {LS}, we can replace $\vp '$ by 
$$\vp =(\n f-1_{Y}¥)\vp '(\n f -¥1_{Y}¥)d(\n f-1_{Y}¥)\vp '(\n f -1_{Y}¥),$$ 
satisfying conditions
(1)--(4).
 
When solving problems in homological or homotopical algebra, one often works
with chain complexes with additional algebraic structure, e.g., chain
algebras or coalgebras.  It is natural to extend the notion of
SDR's to categories of such objects.      

\definition{Definition}  An SDR $(X,d) \sdr{\n}f
(Y,d)\circlearrowleft\vp$ is a {\it 
SDR of chain
(co)algebras} if
\roster
\item $\n $ and $f$ are morphisms of chain (co)algebras, and
\item $\vp $ is a (co)derivation homotopy from $\n f$ to $1_{Y}¥$.
\endroster
\enddefinition

The following notion, introduced by Gugenheim and Munkholm, is 
somewhat weaker than the previous definition for chain coalgebras but perhaps more useful.

\definition{Definition} An SDR $(X,d) \sdr{\n}f (Y,d)\circlearrowleft\vp$ is 
called {\sl  Eilenberg-Zilber (E-Z) data} if $(Y,d, \Delta _{Y}¥)$ and 
$(X,d, \Delta _{X}¥)$ are chain coalgebras and $\n$ is a morphism of 
coalgebras. \enddefinition

Observe that in this case 
$$(d\otimes 1_{Y}¥+1_{Y}¥\otimes d)\bigl((f\otimes f)\Delta _{Y}¥ \vp\bigr)+\bigl((f\otimes 
f)\Delta_{Y}¥ \vp\bigr) 
d=\Delta_{X}¥ f-(f\otimes f)\Delta_{Y}¥,$$
i.e., $f$ is a map 
of coalgebras up to chain homotopy. In fact, $f$ is a DCSH map, as
Gugenheim 
and Munkholm showed  in the following theorem \cite {GM, Thm. 4.1}, which 
proves extremely useful in section 4 of this article.

\proclaim {Theorem 2.1 [GM]}  Let $(X,d) \sdr{\n}f 
(Y,d)\circlearrowleft\vp$ be E-Z data such that $Y$ is simply 
connected and $X$ is connected. Let $F_{1}=f $. Given $F_{i}$ for all $i<k$, let 
$$F_{k}=-\sum _{i+j=k}(F_{i}\otimes F_{j})\Delta _{Y}¥ \vp.$$  Similarly, let $\Phi _{1}= \varphi$, and, given $\Phi _{i}$ for all $i<k$, let 
$$\Phi _{k}=\bigl (\Phi _{k-1}\otimes 1_{Y} + \sum _{i+j=k} \n ^{\otimes 
i}F_{i}\otimes \Phi _{j}\bigr )\Delta _{Y}¥ \vp.$$ Then 
$$\Om(X,d) \sdr{\Om\n}{\tom f} 
\Om(Y,d)\circlearrowleft\tom \vp$$
is an SDR of chain algebras, where $\tom f=\sum _{k\geq 1}(\si)^{\otimes k}F_k s$ and $\tom \vp =\sum _{k\geq 1}(\si)^{\otimes k}\Phi_k s$.
\endproclaim

Let $\cat {EZ}$ be the category with as objects E-Z data $(X,d) \sdr{\n}f 
(Y,d)\circlearrowleft\vp$ such that $Y$ is simply 
connected and $X$ is connected.  A morphism in $\cat {EZ}$  
$$\biggl((X,d) \sdr{\n}f 
(Y,d)\circlearrowleft\vp\biggr)@>>>\biggl ((X',d') \sdr{\n'}{f'} 
(Y',d')\circlearrowleft\vp' \biggr)$$ consists of a pair of morphisms of chain coalgebras $g:(X,d)@>>>(X',d')$ and $h:(Y,d)@>>>(Y',d')$ such that $h\n =\n'g$, $gf=f'h$ and $h\vp=\vp'h$.

\proclaim{Corollary 2.2}  There is a functor $AW: \cat {EZ}@>>> \bigl (\fatcoalg\bigr)^{\rightarrow}$.\endproclaim

\demo{Proof}  Set $AW\bigl((X,d) \sdr{\n}f 
(Y,d)\circlearrowleft\vp\bigr)$  equal to  $\ind ^{-1}(\tom f)\in \fatcoalg (Y,X)$, where $\tom f$ is defined as in Theorem 2.1.  The evident naturality of the definition of $\tom f$ implies that $AW$ is a functor.\qed\enddemo

We call $AW$ the {\sl Alexander-Whitney functor}.

\example{Fundamental Example}  The natural Eilenberg-Zilber and 
Alexander-Whitney equivalences for simplicial
sets provide the most classic example of E-Z data and play a crucial role in the constructions in this article.  Let
$K$ and $L$ be two simplicial sets.  Define morphisms on their normalized
chain complexes
$$\n _{K,L}:C(K)\otimes C(L)@>>>C(K\times L)\quad\text{and}\quad f_{K,L}:C(K\times
L)@>>>C(K)\otimes C(L)$$
by
$$\n _{K,L}(x\otimes y)=\sum \sb {(\mu, \nu)\in \Cal S\sb {p,q}}(-1)^{\sgn (\mu)}(s\sb {\nu \sb
q}...s\sb {\nu \sb 1} x,s\sb {\mu \sb p}...s\sb
{\mu \sb 1} y)$$
where $\Cal S\sb {p,q}$ denotes the set of $(p,q)$-shuffles, $\sgn(\mu)$ is the signature of $\mu$ and  $x\in K\sb
p$, $y\in L\sb q$, and
$$f_{K,L}((x,y))=\sum  \sb {i=0}\sp n \del \sb {i+1}\cdots \del\sb n x\otimes \del
\sb 0\sp i y$$
where $(x,y)\in (K\times L)\sb n$.  We call $\nabla_{K,L} $ the {\it shuffle (or Eilenberg-Zilber) map} and
$f_{K,L}$ the {\it Alexander-Whitney map}.  There is a chain homotopy, $\vp_{K,L}
$, so that  
$$C(K)\otimes C(L)\sdr {\n_{K,L} }{f_{K,L}}C(K\times L)\circlearrowleft\vp_{K,L},
\tag 2.1$$
is an SDR of chain complexes.  Furthermore $\n _{K,L}$ is a map of 
coalgebras, with respect to the usual coproducts, which are defined 
in terms of the natural equivalence $f_{K,L}$.  We have thus defined a functor
$$EZ: \cat {sSet_1}\times\cat{sSet_1}@>>> \cat {EZ},$$
where $\cat{sSet_1}$ is the category of $1$-reduced simplicial sets.  We call $EZ$ the {\sl Eilenberg-Zilber functor}.

When $K$ and $L$ are $1$-reduced, we can apply Theorem 2.1 to the SDR 
(2.1) 
and obtain a new SDR
$$\Om \bigl (C(K)\otimes C(L)\bigr )\sdr {\Om \n _{K,L}}{\tom f_{K,L}}\Om C(K\times 
L)\circlearrowleft\tom \vp_{K,L}.
\tag 2.2$$
In the language of Corollary 2.2, we have a functor from the category $\cat{sSet_1}\times \cat{sSet_1}$  to $\bigl(\fatcoalg\bigr)^\rightarrow$, given by the following composite.
$$\cat{sSet_1}\times \cat{sSet_1}@>EZ>>\cat {EZ}@>AW>> \bigl(\fatcoalg\bigr)^\rightarrow$$

See May's book \cite {M\rm ,
\S 28} and the articles of Eilenberg and MacLane \cite {EM1}, \cite 
{EM2} for further  details.
\endexample 

We can apply our knowledge of this fundamental example to proving the following important result.

\proclaim{Theorem 2.3}  Let $\cat {sSet}_1$ denote the category of $1$-reduced simplicial sets. There is a functor $\widetilde C: \cat {sSet}_1@>>> \cat {wF}$ such that $U\tC=C$, the normalized chains functor, where $U:\cat {wF}@>>>\acoalg$ is the forgetful functor.\endproclaim

\demo{Proof}  Given a $1$-reduced simplicial set $K$,  observe that
$$(AW\circ EZ)(K,K)\in \fatcoalg\bigl(C(K\times K), C(K)\otimes C(K)\bigr).$$  
Define $\Psi_K$ to be the composite
$$\op T\bigl(C(K)\bigr)\acirc \op F@>\op T\bigl((\Delta _K)_\sharp\bigr)\acirc 1>>\op T\bigl(C(K\times K)\bigr)\acirc \op F @>AW\circ EZ(K,K)>> \op T\bigl(C(K)\otimes C(K)\bigr).$$
The pair $\bigl(C(K), \Psi_K\bigr)$ is a weak Alexander-Whitney coalgebra, so we can set
$$\widetilde C(K):=\bigl(C(K), \Psi_K\bigr).$$
It is then immediate that $U\tC(K)=C(K)$.

Given a morphism $h:K@>>>L$ of $1$-reduced simplicial sets, let
$\tC(h)\in \cat F\bigl(\tC(K), \tC(L)\bigr)$ be the morphism of right $\op A$-modules
$$\op T(h_\sharp)\acirc \varepsilon : \op T\bigl(C(K)\bigr)\acirc \op F@>>>\op T\bigl(C(L)\bigr),$$
where $\varepsilon: \op F@>>>\op A$ is the co-unit of $\op F$.  A straightforward diagram chase enables us to establish that $\Psi _L\tC(h)=\bigl(\tC(h)\wedge \tC(h)\bigr)\Psi _K$, ensuring that $\tC(h)$ really is a morphism in $\cat F$.  Key to the success of the diagram chase are the naturality of $AW$ and $EZ$ and of the diagonal map on simplicial sets, as well as the fact that $(\varepsilon \acirc 1)\psi _{\op F}=Id_{\op F}=  (1\acirc \varepsilon )\psi _{\op F} $, i.e., that $\varepsilon$ is a counit for $\psi_{\op F}$.
\qed\enddemo

\head 3. Twisting cochains and twisting functions\endhead

We recall here the algebraic notion of a twisting cochain 
and the simplicial notion of a twisting function, both of which are 
crucial in this article.  We explain the relationship between the
two, which is expressed in terms of a perturbation of the 
Eilenberg-Zilber SDR defined in section 2.  We conclude by recalling 
an important result of Morace and Prout\'e \cite {MP} concerning the 
relationship between Szczarba's twisting cochain and the 
Eilenberg-Zilber equivalence

\definition{Definition} Let $(C,d)$ be a chain coalgebra with 
coproduct $\Delta$, and let 
$(A,d)$ be a chain algebra with product $\mu$.  A {\sl twisting cochain} from $(C,d)$ 
to $(A,d)$ is a degree $-1$ map $t:C@>>>A$ of graded modules such 
that
$$dt+td=\mu (t\otimes t)\Delta.$$
\enddefinition

The definition of a twisting 
cochain $t:C@>>>A$ is formulated precisely so that the following two 
constructions work smoothly. First, let $(A,d)\otimes 
_{t}(C,d)=(A\otimes C, D_{t}),$ where $D_{t}=d\otimes 1_{C}¥+1_{A}¥\otimes d- 
(\mu\otimes 1_{C}¥)(1_{A}¥\otimes t\otimes 1_{C}¥)(1_{A}¥\otimes \Delta)$.  It is easy to 
see that $D_{t}^2=0$, so that $(A,d)\otimes 
_{t}(C,d)$ is a chain complex, which extends $(A,d)$.  

Second, if $C$ 
is connected, let 
$\theta:T\si C_{+}¥@>>>A$ be the algebra map given by $\theta(\si 
c)=t(c)$.  Then $\theta$ is in fact a chain algebra map $\theta :\Om (C,d)@>>>(A,d)$, and  the complex $(A,d)\otimes 
_{t}(C,d)$ is acyclic if and only if $\theta$ is a quasi-isomorphism.  It is equally clear that any algebra map 
$\theta :\Om (C,d)@>>>(A,d)$ gives rise to a twisting cochain via the 
composition
$$C_{+}¥@>\si>>\si C_{+}\hookrightarrow T\si C_{+}@>\theta >> A.$$
In particular, for any two chain coalgebras $(C,d,\Delta)$ and $(C', d', 
\Delta ')$,  the set of DCSH maps from $C$ to $C'$ and the set of twisting cochains from $C$ to $\Om C'$ are naturally in bijective correspondence.

The twisting cochain associated to the cobar construction is a 
fundamental example of this notion.  Let $(C,d,\Delta )$ be a 
simply-connected chain coalgebra.  Consider the linear map
$$t_{\Om C}:C@>>>\Om C:c@>>>\si c.$$
It is a easy exercise to show that $t_{\Om C}$ is a twisting cochain and induces the identity map on $\Om C$.  Thus, in particular, $(\Om C,d)\otimes 
_{t_{\Om}¥}(C,d)$ is acyclic; this is the well-known acyclic cobar 
construction \cite {HMS}.

\definition {Definition/Lemma} Let $t:C@>>>A$ and $t':C'@>>>A'$ be twisting 
cochains.  Let $\ve: C@>>>\Bbb Z$ and $\ve ':C'@>>>\Bbb Z$ be counits, 
and let $\eta :\Bbb Z@>>>A$ and $\eta ':\Bbb Z@>>>A'$ be units. Set 
$$t*t'=t\otimes \eta '\ve '+ \eta\ve\otimes t':C\otimes C'@>>>A\otimes A'$$
Then $t*t'$ is a twisting cochain, called the {\sl cartesian product} 
of $t$ and $t'$.  If $\theta :\Om C@>>>A$ and $\theta ':\Om C'@>>>A'$ 
are the chain algebra maps induced by $t$ and $t'$, then we write 
$\theta*\theta ':\Om (C\otimes C')@>>>A\otimes A'$ for the chain algebra map induced by $t*t'$.
\enddefinition

\remark{Remark}  Observe that Milgram's equivalence $q:\Om (C\otimes C')@>>> \Om C\otimes \Om C'$ is exactly $Id_{\Om C}*Id_{\Om C'}$, which is the chain algebra map induced by $t_{\Om C}*t_{\Om C'}$.\endremark\medskip

\definition {Definition} Let $K$ be a simplicial set and $G$ a simplicial group, 
where the neutral element in any dimension is noted $e$. A 
degree $-1$ map of graded sets $\tau :K@>>>G$ is a {\sl twisting 
function} if 
$$\split
\del _{0}\tau (x)&=\bigl(\tau (\del _{0}x)\bigr)^{-1}\tau (\del 
_{1}x)\\
\del _{i}\tau (x)&=\tau (\del 
_{i+1}x)\quad\text {for all $i>0$}\\
s_{i}\tau (x)&=\tau (s_{i+1}x)\quad \text {for all $i\geq 0$}\\
\tau (s_{0}x)&=e
\endsplit$$
for all $x\in K$.
\enddefinition

The definition of a twisting function $\tau :K@>>>G$ is formulated 
precisely so that if $G$ operates on the left on a simplicial set $L$, 
then we can construct a {\sl twisted cartesian product} of $K$ 
and $L$, denoted $L\times _{\tau}K$, which is a simplicial set such 
that
$(L\times _{\tau}K)_{n}=L_{n}\times K_{n}$, with faces and 
degeneracies given by
$$\split 
\del _{0}(y,x)&=(\tau (x)\cdot \del _{0}y,\del _{0}x)\\
\del _{i}(y,x)&=(\del _{i}y¥,\del 
_{i}x)\quad\text {for all $i>0$}\\
s_{i}(y,x)&=(s _{i}y¥,s_{i}x)\quad \text {for all $i\geq 0$}.
\endsplit$$
If $L$ is a Kan complex, then the projection $L\times _{\tau }K@>>>K$ is 
a Kan fibration \cite {M}.

\example {Example}The canonical twisting 
functions $\lambda _{K}¥ :K@>>>GK:x\mapsto \bar x$ are particularly 
important in this article, in particular because the geometric 
realization of $ GK\times _{\lambda_{K}¥}K$ is acyclic.
\endexample

Twisting cochains and twisting functions are, not surprisingly, very 
closely related.  The theorem below describes their relationship in 
terms of a generalization of the 
Eilenberg-Zilber/Alexander-Whitney equivalences.

\proclaim {Theorem 3.1} For each twisting function $\tau :K@>>>G$ 
there exists a twisting cochain $t(\tau):C(K)@>>>C(G)$ and an SDR
$$C(G)\otimes _{t(\tau)}¥C(K) \sdr{\n_{\tau}¥}{f_{\tau}}¥ C(G\times 
_{\tau}K)\circlearrowleft\vp_{\tau}.$$
Furthermore the choice of $t(\tau)$, $\n_{\tau}$, $f_{\tau}$ and 
$\vp_{\tau}$ can be made naturally.
\endproclaim

Observe that since the realization of $ GK\times _{\lambda_{K}¥}K$ is 
acyclic, $C(GK)\otimes _{t (\lambda_{K})}C(K) $ is acyclic as well, for any natural choice of twisting cochain $t(-)$ fulfilling the conditions of the theorem above. 
Consequently, the induced chain algebra map $\theta (\lambda_{K}): \Om C(K)@>>> C(GK)$ is a quasi-isomorphism.

E. Brown proved the original version of this theorem, for topological 
spaces, by methods of acyclic 
models \cite {B}.  Somewhat later R. Brown \cite {Br} and Gugenheim \cite {G} used 
homological perturbation theory to prove the existence 
of $t(\tau)$ in the simplicial case without defining it explicitly. 
Szczarba was the first to give an explicit, though extremely complex, formula 
for $t (\tau)$, in \cite {Sz}.  

\remark {Convention} Henceforth in this article, the notation 
$sz(\tau)$ will be used exclusively to mean Szczarba's explicit 
twisting cochain, while $sz_{K}$ will always denote $sz(\lambda_{K})$ 
and $Sz _{K}$ the chain algebra map induced by $sz_{K}$.
\endremark\medskip

Recently, in \cite {MP} Morace and Prout\'e provided an alternate, more compact 
construction of $sz(\tau)$, which enabled 
them to prove that  $sz_{K}$ commutes with the shuffle map, 
as described below.

\proclaim {Theorem 3.2 [MP]} Let $K$ and $L$ be reduced simplicial 
sets.  Let $\rho : G(K\times L)@>>>GK\times GL$ denote the 
homomorphism of simplicial groups defined in the introduction. Then the diagram of graded module 
maps 
$$\xymatrix
{C(K)\otimes C(L) \ar[dd] ^{sz_{K}*sz_{L}}\ar[rr]^{\n 
_{K,L}}&&C(K\times L)\ar [d]^{sz_{K\times L}}\\
&&C(G(K\times L))\ar [d]^{\rho _{\sharp}}\\
C(GK)\otimes C(GL)\ar [rr]^{\n _{GK,GL}}&&C(GK\times GL)}$$
commutes.
\endproclaim

The following 
corollary of Theorem 3.2 is crucial to the development in the next section.  

\proclaim {Corollary 3.3} Let $K$ and $L$ be $1$-reduced simplicial 
sets, and let $\rho$ be as above.  Then the diagram of chain algebra 
maps
$$\xymatrix
{\Om C(K\times L)\ar [rr]^{\tom f_{K,L}}\ar [d]^{Sz_{K\times L}}&&
\Om \bigl(C(K)\otimes C(L)\bigr )\ar[dd] ^{Sz_{K}*Sz_{L}}\\
C(G(K\times L))\ar [d]^{\rho _{\sharp}}&&\\
C(GK\times GL)\ar [rr]^{f_{GK,GL}}&&C(GK)\otimes C(GL)}$$

commutes up to homotopy of chain algebras.
\endproclaim

\demo {Proof} Recall 
that $\n$ is always a map of coalgebras, so that it induces a map of 
chain algebras $\Om \n$ on cobar constructions.
As an immediate consequence of Theorem 3.2, we obtain that the 
diagram of chain algebras
$$\xymatrix
{\Om \bigl(C(K)\otimes C(L)\bigr ) \ar[dd] ^{Sz_{K}*Sz_{L}}
\ar[rr]^{\Om \n _{K,L}}&&\Om C(K\times L)\ar [d]^{Sz_{K\times L}}\\
&&C(G(K\times L))\ar [d]^{\rho _{\sharp}}\\
C(GK)\otimes C(GL)\ar [rr]^{\n _{GK,GL}}&&C(GK\times GL)}$$
commutes.  It suffices to check the commutativity for generators of 
$\Om \bigl(C(K)\otimes C(L)\bigr )$, i.e, for elements of $\si 
\bigl(C(K)\otimes C(L)\bigr )_{+}$, which is equivalent to the 
commutativity of the diagram in Theorem 3.2.

If $\Phi =f_{GK,GL}\circ\rho_{\sharp}\circ Sz_{K\times L}\circ\tom \vp _{K,L}:\Om 
C(K\times L)@>>>C(GK)\otimes C(GL)$, then
$$\split
(d\otimes 1+1\otimes d)\Phi +\Phi 
d&=f_{GK,GL}\circ\rho_{\sharp}\circ Sz_{K\times 
L}\circ\Om\n\circ\tom f_{K,L}-f_{GK,GL}\circ\rho_{\sharp}\circ Sz_{K\times L}\\
&=f_{GK,GL}\circ\n _{GK,GL}\circ Sz_{K}* Sz_{L}\circ\tom f_{K,L}
-f_{GK,GL}\circ\rho_{\sharp}\circ Sz_{K\times L}\\¥
&=Sz_{K}* Sz_{L}\circ\tom f_{K,L}
-f_{GK,GL}\circ\rho_{\sharp}\circ Sz_{K\times L}.
\endsplit$$
The map $\Phi$ is thus a chain homotopy from $Sz_{K}*Sz_{L}\circ\tom 
f_{K,L}$ to $f_{GK,GL}\circ\rho_{\sharp}\circ Sz_{K\times L}$.  
Furthermore, since, according to Theorem 2.1,  $\tom \vp$ is a $(\Om \n\circ \tom f, 
1)$-derivation, $\Phi$ is a $(Sz_{K}*Sz_{L}\circ\tom f_{K,L}, 
f_{GK,GL}\circ\rho_{\sharp}\circ Sz_{K\times 
L})$-derivation.  The diagram in the statement of the theorem commutes 
therefore up to homotopy of chain algebras.
\qed\enddemo


\head 4. The canonical Adams-Hilton model\endhead

Our goal in this section is to define and establish the key properties 
of the Alexander-Whitney cobar diagonal.
Throughout the section we abuse notation slightly and write $f_{K}$, $\n _{K}$ and $\vp_{K}$ 
instead of $f_{K,K}$, $\n _{K,K}$ and $\vp_{K,K}$.  Recall furthermore the functors $\tom: \cat{wF}@>>>\cat {wH}$ (Lemma 1.4) and $\tC:\cat {sSet}_1@>>>\cat {wF}$ (Theorem 2.3).

\definition {Definition} Let $K$ be a $1$-reduced simplicial set.  The {\sl canonical Adams-Hilton model} for $K$ is $\tom\tC (K)$.  The coproduct $\psi_K$ on the canonical Adams-Hilton model is called the {\sl Alexander-Whitney (A-W) cobar diagonal}.\enddefinition

Unrolling the definition of $\psi_K$, we see that it is equal to the following composite.
$$\Om C(K)@>\Om (\Delta _{K})_{\sharp}>>\Om C(K\times K)@>\tom 
f_{K}>>\Om (C(K)\otimes C(K))@>q>>\Om C(K)\otimes \Om C(K).$$
We show in the next 
two results that it is, in particular, cocommutative up to 
homotopy of chain algebras and strictly coassociative.  Thus, $\tom\tC(K)\in \cat H$, i.e., $\tC(K)$ is a {\sl strict} Alexander-Whitney coalgebra. 

\proclaim {Proposition 4.1}The Alexander-Whitney cobar diagonal $\psi_{K}$ is 
cocommutative up to homotopy of chain algebras for all $1$-reduced simplicial sets $K$.
\endproclaim

\demo {Proof}
Consider the following diagram, in which $sw$ denotes both the 
simplicial coordinate switch map and the algebraic tensor switch map.
$$\xymatrix
{\Om C(K) \ar [dr]_{\Om (\Delta)_{\sharp}}\ar [r]^{\Om 
(\Delta)_{\sharp}}&\Om C(K\times K)\ar [d]^{\Om (sw)_{\sharp}}\ar 
[r]^{\tom f}&\Om (C(K)\otimes C(K))\ar[d]^{\Om (sw)}\ar [r]^q&\Om 
C(K)\otimes \Om C(K)\ar [d]^{sw}\\
&\Om C(K\times K)\ar [r]^{\tom f}&\Om (C(K)\otimes C(K))\ar [r] ^q&\Om 
C(K)\otimes \Om C(K)}$$
The triangle on the left and the square on the right commute for obvious 
reasons, while the middle square commutes up to chain homotopy, as
$$\split 
\Om (sw)\circ \tom f=&\tom f\circ \Om \n \circ \Om (sw)\circ \tom f\\
=&\tom f\circ \Om (sw)_{\sharp}\circ \Om \n\circ  \tom f\qquad\text {since $\n \circ sw= 
(sw)_{\sharp}\circ \n $} \\
\simeq&\tom f\circ \Om (sw)_{\sharp}.
\endsplit$$
The homotopy in the last step is provided by $\tom f\circ \Om 
(sw)_{\sharp}\circ \tom \vp$. Hence, the whole diagram commutes up to chain 
homotopy, where $q\circ \tom f\circ \Om (sw)_{\sharp}\circ \tom 
\vp\circ \Om \Delta _{\sharp}$ provides the necessary homotopy.
\qed\enddemo

\proclaim {Theorem 4.2} The Alexander-Whitney cobar diagonal $\psi_{K}$ is 
strictly coassociative for all $1$-reduced simplicial sets $K$. 
\endproclaim

\demo{Proof} We need to show that $(\psi _{K}\otimes 1)\psi 
_{K}=(1\otimes \psi _{K})\psi _{K}$, which means that we need to show 
that the following diagram commutes. (Note that we drop the subscript 
$K$ for the remainder of this proof.)
$$\xymatrix
{\Om C(K)\ar[d]^{\Om (\Delta )_{\sharp}}\ar[rr]^{\Om (\Delta 
)_{\sharp}}&&
\Om C(K^2)\ar [rr]^{\tom f}&&
\Om (C(K)^{\otimes 2})\ar [r]^q&
(\Om C(K))^{\otimes 2}\ar [d]^{1\otimes \Om (\Delta)_{\sharp}}\\
\Om C(K^2)\ar [d]^{\tom f}&&&&&
\Om C(K)\otimes \Om C(K^2)\ar [d]^{1\otimes {\tom f}}\\
\Om (C(K)^{\otimes 2})\ar [d]^q&&&&&
\Om C(K)\otimes (\Om C(K)^{\otimes 2})\ar [d]^{1\otimes q}\\
(\Om C(K))^{\otimes 2}\ar [rr]^{\Om (\Delta )_{\sharp}\otimes 1}&&
\Om C(K^2)\otimes \Om C(K)\ar [rr]^{\tom f\otimes 1}&&
\Om (C(K)^{\otimes 2})\otimes \Om C(K)\ar [r]^{q\otimes 1}&
(\Om C(K))^{\otimes 3}
}$$
In order to prove that the square above commutes, we divide it into nine 
smaller squares
$$\xymatrix
{\Om C(K)\ar[d]^{\Om (\Delta )_{\sharp}}\ar[rr]^{\Om (\Delta )_{\sharp}}&&
\Om C(K^2)\ar [d] ^{\Om (1\times \Delta)_{\sharp}}¥\ar [rr]^{\tom f}&&
\Om (C(K)^{\otimes 2})\ar [d]^{\Om (1\otimes (\Delta )_{\sharp})}\ar [r]^q&
(\Om C(K))^{\otimes 2}\ar [d]^{1\otimes \Om (\Delta)_{\sharp}}\\
\Om C(K^2)\ar [d]^{\tom f}\ar [rr] ^{\Om (\Delta\times 1)_{\sharp}}&&
\Om C(K^3)\ar [d]^{\tom f_{K^2,K}}\ar [rr]^{\tom f_{K,K^2}}&&
\Om (C(K)\otimes C(K^2))\ar [d]^{\tom (1\otimes f)}\ar [r]^q&
\Om C(K)\otimes \Om C(K^2)\ar [d]^{1\otimes {\tom f}}\\
\Om (C(K)^{\otimes 2})\ar [rr]^{\Om ((\Delta)_{\sharp}\otimes 1)}\ar [d]^q&&
\Om (C(K^2)\otimes C(K))\ar [d]^q\ar [rr]^{\tom (f\otimes 1)}&&
\Om (C(K)^{\otimes 3})\ar [d]^q\ar [r]^q&
\Om C(K)\otimes (\Om C(K)^{\otimes 2})\ar [d]^{1\otimes q}\\
(\Om C(K))^{\otimes 2}\ar [rr]^{\Om (\Delta )_{\sharp}\otimes 1}&&
\Om C(K^2)\otimes \Om C(K)\ar [rr]^{\tom f\otimes 1}&&
\Om (C(K)^{\otimes 2})\otimes \Om C(K)\ar [r]^{q\otimes 1}&
(\Om C(K))^{\otimes 3}
}$$
and show that each of the small squares commutes, with one exception, 
for which we can correct.  We label each 
small square with its row and column number, so that, e.g., square 
$(2,3)$ is 
$$\xymatrix
{\Om (C(K)\otimes C(K^2))\ar [d]^{\tom (1\otimes f)}\ar [r]^q&
\Om C(K)\otimes \Om C(K^2)\ar [d]^{1\otimes {\tom f}}\\
\Om (C(K)^{\otimes 3})\ar [r]^q&
\Om C(K)\otimes (\Om C(K)^{\otimes 2}).}$$

The commutativity of eight of the nine small squares is immediate. Square $(1,1)$ commutes since 
$\Delta$ is coassociative.  Squares 
$(1,2)$ and $(2,1)$ commute by naturality of $f$, while squares 
$(1,3)$ and $(3,1)$ commute by naturality of $q$.  The commutativity of squares $(2,3)$ and $(3,2)$ is an immediate consequence of Proposition 2.2. Finally, a simple 
calculation shows that square $(3,3)$ commutes as well.

Let $q^{(2)}=(1\otimes q)q=(q\otimes 1)q$.  In the case of square $(2,2)$, we show that 
$$Im \bigl(\tom (1\otimes f)\circ \tom f_{K,K^2}-\tom (f\otimes 1)\circ 
\tom f_{K^2,K}\bigr )\subseteq \ker q^{(2)},$$
which suffices to conclude that the large square commutes, since we 
know that the other eight small squares commute. 

Let $c_{1,2}$ and 
$c_{2,1}$ denote the usual coproducts on $C(K)\otimes C(K^2)$ and 
$C(K^2)\otimes C(K)$, respectively. Given any 
$z\in C(K^3)$, use the Einstein summation convention in writing $f _{K,K^2}(z)=x_{i}\otimes y^i$, 
$c_{K}(x_{i})=x_{i,j}\otimes x_{i}^j$, and $c_{K^2}(\vp (y^i))=\vp 
(y^i)_{k}\otimes \vp (y^i)^k$, so that
$$\align
q^{(2)}\bigl(\si (1\otimes f)\bigr )^{\otimes 2} 
c_{1,2}&(1\otimes \vp)f _{K,K^2}(z)\\
&=q^{(2)}\bigl(\si (1\otimes f)\bigr )^{\otimes 2} 
c_{1,2}\bigl (x_{i}\otimes \vp (y^i)\bigr )\\
&=q^{(2)} \bigl(\si (1\otimes f)\bigr )^{\otimes 
2}\bigl (\pm x_{i,j}\otimes \vp (y^i)_{k}\otimes x_{i}^j\otimes \vp 
(y^i)^k\bigr )\\
&=q^{(2)}\biggl (\pm \si \bigl (x_{i,j}\otimes f (\vp (y^i)_{k})\bigr 
)\si \bigl (x_{i}^j 
\otimes f(\vp (y^i)^k)\bigr )\biggr)\\
&=(1\otimes q)\bigl((\si x_i\otimes 1)(1\otimes \si f\vp(y^i))\pm (1\otimes \si f\vp (y^i))(\si x_i\otimes 1)\bigr)
\endalign$$
since $(1\otimes q)(\si (u\otimes v))=0$ unless $|u|=0$ or $|v|=0$.  
This last sum is $0$, however, since $f\vp=0$.

Similarly, $q^{(2)}\bigl(\si (f\otimes 1)\bigr )^{\otimes 2} 
c_{2,1}(\vp\otimes 1)f _{K^2,K}(z)=0$.  Applying Gugenheim and 
Munkholm's formula from Theorem 2.1, we obtain for all $z\in C(K^3)$
$$\align 
q^{(2)}\tom (1\otimes f)\tom f_{K,K^2}(\si z)&=q^{(2)}\si(1\otimes 
f)f_{K,K^2}(z)\\
&=q^{(2)}\si (f\otimes 1)f_{K^2,K}(z)\\
&=q^{(2)}\tom (f\otimes 1) 
\tom f_{K^2,K}(\si z),\endalign$$
since in general 
$$(1\otimes f_{L,M})f_{K,L\times M}=(f_{K,L}\otimes 
1)f_{K\times L, M}:C(K\times L\times M)@>>>C(K)\otimes C(L)\otimes 
C(M).\qed$$
\enddemo

\proclaim {Proposition 4.3}The chain algebra quasi-isomorphism $Sz_{K}:\Om 
C(K)@>>>C(GK)$ induced by Szczarba's twisting cochain $sz_{K}$ is a 
map of chain coalgebras up to homotopy of chain algebras, with 
respect to the Alexander-Whitney cobar diagonal and the usual coproduct 
$c_{GK}=¥f_{GK}\circ (\Delta _{GK})_{\sharp}$ on 
$C(GK)$, i.e., the diagram
$$\xymatrix
{\Om C(K)\ar [rr]^{\psi _{K}}\ar [d]^{Sz_{K}}&&\Om C(K)\otimes 
\Om C(K)\ar [d] ^{Sz_{K}\otimes Sz_{K}}\\
C(GK)\ar [rr]^{c_{GK}¥}&&C(GK)\otimes C(GK)}$$
commutes up to homotopy of chain algebras.
\endproclaim

\remark{Remark}  Since $c_{GK}$ is homotopy cocommutative, Proposition 4.3 
implies immediately that $\psi _{K}$ is homotopy cocommutative as well.  
We consider, however, that it is worthwhile to establish the homotopy 
cocommutativity of $\psi _{K}$ independently, as we do in Proposition 4.1, 
since we obtain an explicit 
formula for the chain homotopy.
\endremark\medskip

\demo {Proof} We can expand and complete the diagram in the statement of the 
theorem to obtain the diagram below.
$$\xymatrix
{\Om C(K)\ar [d]^{\Om (\Delta _{K})_{\sharp}}\ar [rrr]^{Sz_{K}}&&&
C(GK)\ar [d]^{(G\Delta _{K})_{\sharp}}\ar@/^5pc/ 
[dd] ^{(\Delta _{GK})_{\sharp}}   \\
\Om C(K\times K)\ar [d] ^{\tom f_{K}}\ar [rrr]^{Sz_{K^2}}&&&C\bigl (G(K\times K)\bigr)\ar [d]^{\rho _{\sharp}}\\
\Om (CK\otimes CK)\ar [d] ^q\ar@<1ex> [drrr]^{
{Sz_{K}*Sz_{K}}}&&&C(GK\times GK)\ar [d] ^{f_{GK}}\\
\Om CK\otimes \Om CK\ar [rrr]^{Sz_{K}\otimes Sz_{K}}&&&C(GK)\otimes C(GK)}$$¥
The top square commutes exactly, by naturality of the twisting 
cochains.  An easy calculation shows that the bottom triangle commutes 
exactly.  Since Corollary 3.3 implies that the middle square commutes up 
to homotopy of chain algebras, we can conclude that the theorem is 
true.  In particular  $f_{GK}\rho_{\sharp}Sz_{K\times K}\tom \vp _{K}
\Om (\Delta _{K})_{\sharp}$ is an appropriate derivation homotopy.
\qed\enddemo

It would be interesting to determine under what conditions $Sz_{K}$ is a strict map of Hopf algebras.  We have checked 
that $(Sz_{K}\otimes Sz_{K})\psi _{K}=c_{GK}Sz_{K}$ up through degree $3$ and will show in a later paper \cite {HPS2} that $Sz _K$ is a strict Hopf algebra map when $K$ is a suspension.

Even if $Sz_{K}$ is not a strict coalgebra map, we know that 
it is at least the next best thing, as stated in the following theorem.

\proclaim {Theorem 4.4} Any natural map $\theta_{K}:\Om 
C(K)@>>> C(GK)$ of chain algebras is a DCSH map, with respect to any natural choice of strictly coassociative coproduct $\chi _K$ on $\Om C(K)$.\endproclaim

\remark{Remark}  Proposition 4.3 is, of course, an immediate corollary of 
Theorem 4.4.  The independent proof of Proposition 4.3 serves to provide 
an explicit formula for the homotopy between $Sz_{K}\otimes 
Sz_{K})\psi _{K}$ and $c_{GK}Sz_{K}$.  The proof below sacrifices all hope of explicit formulae on the altar of extreme generality. 
\endremark\medskip

\demo{Proof} Let $\overline \Delta [n]$ denote the quotient of the 
standard simplicial $n$-simplex $\Delta [n]$ by its 0-skeleton. Recall 
from \cite {MP} that there is a contracting chain homotopy 
$\bar h:C(G\overline \Delta [n])@>>>C(G\overline \Delta [n])$.   The functor $C\bigl(G(-)\bigr)$ from reduced simplicial sets to connected chain algebras is therefore acyclic on the set of models $\frak M=\{\overline \Delta [n]\mid n\geq 0\}$.

On the other hand, the functor $C$ from reduced simplicial sets to connected chain coalgebras is free on $\frak M$.  In particular,  the set $\{ \iota _n\in C(\overline \Delta [n])\mid n\geq 0\}$ gives rise to basis of $C(K)$ for all $K$, where $\iota_n$ denotes the unique nondegenerate $n$-simplex of $\overline \Delta 
[n]$. 

Theorem 4.4 therefore follows immediately from Theorem 1.5, where $\bold D$ is the category of reduced simplicial sets, $X= (\Om C(-), \chi_{-})$ and $Y=C\bigl(G(-)\bigr)$.\qed\enddemo

\remark{Remark}  The results in this section beg the question of iteration of the cobar construction.  Let $K$ be any $2$-reduced simplicial set.  Since $\psi_K$ is strictly coassociative and $Sz_K$ is a DCSH map, we can apply the cobar construction to the quasi-isomorphism $Sz_K:\tom \tC(K)=(\Om C(K),\psi _K)@>>>C(GK)$ and consider the composite
$$\Om\tom\tC(K)@>\tom Sz_K>>\Om C(GK)@>Sz_{GK}>>C(G^2K),$$
which is a quasi-isomorphism of chain algebras (see Lemma 1.4 and Theorem 2.3 for explanation of the notation $\tom $ and $\tC$).
The question is now whether there is a canonical, topologically-meaningful way to define a coassociative coproduct on $\Om\tom\tC(K)$, in order to iterate the process. In other words, is there a natural, coassociative coproduct on $\Om\tom\tC(K)$ with respect to which $Sz_{GK}\circ \tom Sz_K$ is a DCSH map? Equivalently, does $\tom\tC(K)$ admit a natural Alexander-Whitney coalgebra structure, with respect to which $Sz_K$ is a morphism in $\cat F$? 

Using the notion of the {\sl diffraction} from \cite{HPS} and the more general version of the Cobar Duality Theorem proved there, we can show that $\tom\tC(K)$ admits a natural {\sl weak} Alexander-Whitney structure, with respect to which $Sz_K$ is a morphism in $\cat {wF}$.  Consequently, $\Om\tom\tC(K)$ indeed admits a natural coproduct, but it is not necessarily coassociative, which prevents us from applying the cobar construction again.

In \cite{HPS2} we show that if $EK$ is the suspension of a  $1$-reduced simplicial set $K$, then $\tom\tC(EK)$ does admit a natural, strict Alexander-Whitney coalgebra structure and that $Sz_{EK}$ is then a morphism in $\cat F$.  We conjecture that this result generalizes to higher suspensions and correspondingly higher iterations of the cobar construction. 
\endremark

\head 5. The Baues coproduct and the Alexander-Whitney cobar diagonal\endhead

We show in this section that the coproduct defined by Baues in \cite {Ba} is the same as the Alexander-Whitney cobar diagonal defined in section 4.  As mentioned in the introduction, this result is at first sight quite surprising, since there is an obvious asymmetry in Baues's combinatorial definition, which is well hidden in our definition of the A-W cobar diagonal.

We begin by recalling the definition of Baues's coproduct on $\Om C(K)$, where $K$ is a $1$-reduced simplicial set.  For any $m\leq n\in \Bbb N$, let $[m,n]=\{j\in \Bbb N\mid m\leq j\leq n\}$. Let $\bold \Delta$ denote the category with objects
$$Ob \bold \Delta=\{[0,n]\mid n\geq 0\}$$
and 
$$\bold \Delta \bigl ([0,m], [0,n]\bigr)=\{ f:[0,m]@>>>[0,n]| f\text{ order-preserving set map}\}.$$
Viewing the simplicial set $K$ as a contravariant functor from $\bold \Delta $ to the category of sets, given $x\in K_n:=K([0,n])$ and $0\leq a_1<a_2<\cdots <a_m\leq n$, let
$$x_{a_1...a_m}:=K(\bold a)(x)\in K_m$$
where $\bold a:[0,m]@>>> [0,n]:j\mapsto a_j$.

Let $x$ be a nondegenerate $n$-simplex of $K$. Baues's coproduct $\widetilde \psi$ on $\Om C(K)$ is defined by
$$\widetilde\psi (\si x)=\sum \Sb 0\leq m<n\\0< a_1<\cdots<a_m<n\endSb (-1)^{\ell (\bold a)}\si x_{0...a_1}\si x_{a_1...a_2}\cdots\si x_{a_m...n}\otimes \si x_{0a_1...a_m n}$$
where
$$\ell (\bold a)=(a_1-1)+\bigl(\sum _{i=2}^m (i-1)(a_i-a_{i-1}-1)\bigr) +m(n-a_m-1).$$
Baues showed in \cite {Ba} that $\widetilde\psi$ was strictly coassociative and that it was cocommutative up to derivation homotopy, providing an explicit derivation homotopy for the cocommutativity.

To prove that the A-W cobar diagonal agrees with the Baues coproduct, we examine closely each summand of $\tom f_{K,L}$, determining precisely what survives upon composition with $q:\Om (CK\otimes CL)@>>>\Om CK\otimes \Om CL$. Henceforth, in the interest of simplifying the notation, we no longer make signs explicit.

Recall from section 2 the Eilenberg--Zilber SDR of normalized chain complexes
$$C(K)\otimes C(L)\sdr {\n_{K,L} }{f_{K,L}}C(K\times L)\circlearrowleft\vp_{K,L},
$$
where we can rewrite the definitions of $f_{K,L}$ and $\n _{K,L}$ in terms the notation introduced above as follows. In degree $n$, 
$$\align
(\n_{K,L}) _n(x\otimes y)&=\sum _{\ell =0}^n\sum \Sb A\cup B=[0, n-1]\\|A|=n-\ell, |B|=\ell\endSb \pm (s_A x, s_By)\tag 5.1\\
\intertext{and}
(f_{K,L})_n(x,y)&=\sum _{\ell =0}^nx_{0...\ell}\otimes y_{\ell...n}\tag 5.2
\endalign$$
where $s_I$ denotes $s_{i_r}\cdots
s_{i_2}s_{i_1}$ for any set $I$ of non-negative integers
$i_1<i_2<\cdots <i_r$ and $| I|$ denotes the cardinality of $I$. 
There is also a recursive formula for $\varphi_{K,L}$, due to Eilenberg and MacLane \cite {EM2}. Let $g=\n_{K,L} f_{K,L}$. Then
$$(\varphi _{K,L})_n=-(g)'s_0+\bigl((\varphi _{K,L})_{n-1}\bigr )',\tag 5.3$$
where the prime denotes the derivation operation on simplicial operators, i.e.,
$$h=s_{j_n}\cdots s_{j_0}\del _{i_0}\cdots \del _{i_m}\Rightarrow h'=s_{j_n+1}\cdots s_{j_0+1}\del _{i_0+1}\cdots \del _{i_m+1}.$$

Let $\phihat$ denote the degree +1 map 
$$\phihat:
C(K\times L)@>\varphi_{K,L}>>C(K\times L)@>(\Delta _{K\times L})_{\sharp}>> 
C((K\times L)^2)@>f_{K\times L, K\times L}>>
C(K\times L)^{\otimes 2}.$$
If $K$ and $L$ are 0-reduced, consider the pushout $CK\vee CL$ of
the complexes $CK$ and
$CL$ over $\Bbb Z$ and the map of chain complexes
$$\kappa:CK\otimes CL\longrightarrow CK\vee CL$$
defined by $\kappa(x\otimes 1)=x$, $\kappa(1\otimes y)=y$ and  $\kappa(x\otimes
y)=0$ if $|x|,|y|>0$. Define a family of linear maps
$$\overline{\Cal F}=\{{\Fbar}_k:C(K\times L)\longrightarrow (CK\vee CL)^{\otimes k}\mid \deg \Fbar =k-1, k\geq 1\}$$ 
by 
$$\align
&{\Fbar}_1:
C(K\times L)@>f_{K,L}>>CK\otimes CL@>\kappa>> CK\vee CL\tag 5.4 \\ 
&{\Fbar}_k:C(K\times L)@>\phihat>>C(K\times L)^{\otimes 2}
@>{-\sum _{i+j=k}\Fbar _i\otimes \Fbar _j}>>(CK\vee CL)^{\otimes k}.\tag 5.5
\endalign$$
We can use the family $\overline {\Cal F}$ to obtain a useful factorization of $q\circ \tom f_{K,L}$, as follows.

Observe first, by comparison with the construction given after the statement of Theorem 2.1, that 
$$\bar F_k=\kappa ^{\otimes k}F_k\tag 5.6$$
where $F_k$ is defined as in Theorem 2.1. Next note that for any pair of $1$-connected chain coalgebras $(C,d)$ and $(C',d')$, the algebra map $\gamma :T\si (C_+\vee C_+')@>>>T\si C_+\otimes T\si C'_+$ specified by $\gamma (\si c)=\si c\otimes 1$ and $\gamma (\si c')=1\otimes \si c'$ for $c\in C$ and $c'\in C'$ commutes with the cobar differentials, i.e., it is a chain algebra map $\gamma:\Om ((C,d)\vee(C',d'))@>>>\Om (C,d)\otimes \Om (C',d')$.  Furthermore,
$$q=\gamma\circ T(\si \kappa s):\Om (CK\otimes CL)@>>>\Om CK\otimes \Om CL,$$
which implies, by (5.6), that when $K$ and $L$ are $1$-reduced, 
$$q\tom f_{K,L}=\gamma \circ \sum _{k\geq 1}(\si )^{\otimes k}\bar F_ks:\Om C(K\times L)@>>>\Om C(K)\otimes C(L).$$ 
In particular, the A-W cobar diagonal is equal to $ \gamma\circ \sum _{k\geq 1}(\si )^{\otimes k}\bar F_ks\circ \Om (\Delta _K)_{\sharp}$.

Thanks to the decomposition above, we obtain as an immediate consequence of the next theorem that 
$$\psi _K =\widetilde\psi.\tag 5.7$$

\proclaim{Theorem 5.1}
Let $\overline{\Cal F}=\{{\Fbar}_k:C(K\times L)\longrightarrow (CK\vee CL)^{\otimes k}\mid \deg \Fbar =k-1, k\geq 1\}$ denote the family defined above. If $j\ge2$ then
$$({\Fbar}_i\otimes {\Fbar}_j)\,\phihat=0\tag 5.8$$
If $K$ and $L$ are 1-reduced then
$${\Fbar}_k(x,y)=
\sum_{\{0<i_1<\cdots<i_{r}<n\}} 
\pm y_{0i_1\ldots i_{r} n}\otimes
x_{0\ldots i_1}\otimes\cdots\otimes x_{i_{r}\ldots n}\tag 5.9$$
In this summation we adopt the convention that
1-simplices $x_J$ or $y_J$ for $|J|=2$ are to be identified to the
unit (\sl {not} to zero), and
consider only those terms for which exactly $k$ non-trivial tensor factors remain.
\endproclaim

To see how Theorem 5.1 implies (5.7), note that (5.9) implies that for all $x\in C_nK$,
$$\multline
\gamma \biggl(\sum _{k\geq 1}(\si )^{\otimes k}\Fbar _k(x,x)\biggr)\\= \sum_{\{0<i_1<\cdots<i_{r}<n\}} 
\pm 
\si x_{0\ldots i_1}\cdots \si x_{i_{r}\ldots n}\otimes \si x_{0i_1\ldots i_{r} n}.\endmultline$$

In the proof of Theorem 5.1 we rely on the following lemmas.

The map ${\Fbar}_1$ is easy to identify, by definition of $f$ and $\kappa $.
\proclaim{Lemma 5.2}
For all $(x,y)$ in $K\times L$, we have ${\Fbar}_1(x,y)=x+y$.
\endproclaim

To understand ${\Fbar}_k$ for $k\ge2$, we use the following explicit
formula for $\varphi$.

\proclaim{Lemma 5.3}
Let $A,B$ be disjoint sets such that $ A\cup B=[m+1,n]$ and $|B|=r-m$ for $r>m$.
For $(x,y)\in K_n\times L_n$, let
$$\varphi^{A,B}(x,y)=(s_{A\cup\{m\}}\,x_{0\ldots r},\,s_B\,y_{0\ldots mr\ldots n})
\in K_{n+1}\times L_{n+1}.$$
Then the Eilenberg--Zilber homotopy $\varphi$ is given by
$$\varphi(x,y) =
\sum\Sb m<r\\ A\cup B=[m+1,n] \\
    |A|=n-r,\;|B|=r-m\endSb 
\pm\varphi^{A,B}(x,y).$$
\endproclaim


The following result should be compared with Lemma 5.2 and the classical result $\varphi^2=0$.

\proclaim{Lemma 5.4}
For a general term $\varphi^{A,B}(x,y)$ of Lemma 5.3 and
$0\le\ell\le n+1$,
$$\align
\Fbar_1((\varphi^{A,B}(x,y))_{0\ldots\ell})&=
\left\{\alignedat 2
&x_{0\ldots\ell}\;\;+\;\;y_{0\ldots\ell}
&&:\ell\le m,\\ 
&y_{0\ldots m r\ldots r-m+\ell-1}&&
:\ell>m \text{ and } [m+1,\ell-1]\subseteq A,\\ 
&0&&:\text{else} 
\endalignedat \right.
\\
\Fbar_1((\varphi^{A,B}(x,y))_{\ell\ldots n+1})
&=\left\{\alignedat 2
&x_{r-n+\ell-1\ldots r}&&:\ell>m\text{ and } [\ell,n]\subseteq B,\\  
&y_{\ell-1\ldots n}&&:\ell>m\text{ and } [\ell,n]\subseteq A,\\ 
&0&&:\text{else.} 
\endalignedat\right.
\\
\varphi((\varphi^{A,B}(x,y))_{0\ldots\ell})
&=\quad 0\qquad:\ell>m \text{ and } [m+1,\ell-1]\not\subseteq A
\\
\varphi((\varphi^{A,B}(x,y))_{\ell\ldots n+1})
&=\quad 0\qquad \text{ always.}
\endalign$$
\endproclaim

\demo{Proof of Theorem 5.1}
If $j\ge2$, then we have
$$
\Fbar_j((\varphi^{A,B}(x,y))_{\ell\ldots n+1})
\;\;=\;\;
-\sum_{j_1+j_2=j}(\Fbar_{j_1}\otimes\Fbar_{j_2})\,f\,\Delta _\sharp\,
\varphi((\varphi^{A,B}(x,y))_{\ell\ldots n+1})\;\;=\;\;0
$$
by the final result of Lemma 5.4, and so
$$(\Fbar_i\otimes\Fbar_j)\phihat(x,y)\;\;=\;\;
\sum_{\ell,m,r,A,B}\!\!\!
\Fbar_i((\varphi^{A,B}(x,y))_{0\ldots\ell})
\;\otimes\;
\Fbar_j((\varphi^{A,B}(x,y))_{\ell\ldots n+1})
\;\;=\;\;0,$$
proving the first part of the theorem. 

For the second part, note that for $k=1$ the right hand side
of (5.9) reduces to 
$$y_{0n}\otimes x_{0\ldots n}
\;\;+\;\;
y_{0 \ldots n}\otimes x_{0 1}\otimes\cdots\otimes x_{n-1 n} $$
which we identify with $x+y=\Fbar_1(x,y)$.

For $k=2$ we use the first two results of Lemma 5.4 and the fact
that $B\neq\emptyset$ so $A$ cannot contain both $[m+1,\ell-1]$ and
$[\ell,n]$, to establish that 
$$\align
\Fbar_2(x,y)
&=\sum_{\ell,m,r,A,B}\!\!\!
\Fbar_1((\varphi^{A,B}(x,y))_{0\ldots\ell})
\otimes
\Fbar_1((\varphi^{A,B}(x,y))_{\ell\ldots n+1})\\ 
&=\sum_{\ell,m,r}\quad
y_{0\ldots mr\ldots r-m+\ell-1}
\,\otimes\,
x_{r-n+\ell-1\ldots r}\,.\endalign$$
Here $[m+1,\ell-1]=A$ and $[\ell,n]=B$, so $r-m=|B|=n-\ell+1$ and we
have
$$\Fbar_2(x,y)\;\;=\;\;
\sum_{m<r}\;
y_{0\ldots m r\ldots n}
\;\otimes\;
x_{m \ldots r}$$
which agrees with (5.9). For $k\geq 2$ we have, by (5.8),
$$\Fbar_{k+1}(x,y)
\;\;=\sum\Sb \ell,m,r,A,B\\ i+j=k\endSb\!\!\!
\left((\Fbar_{i}\otimes\Fbar_{j})\,f\Delta_\sharp\,\varphi
((\varphi^{A,B}(x,y))_{0 \ldots \ell})\right)
\;\otimes\;
\Fbar_1((\varphi^{A,B}(x,y))_{\ell \ldots n+1})$$
and Lemma 5.4 tells us once again we must take
$A=[m+1,\ell-1]$, $B=[\ell,n]$, $\ell=n+m-r+1$ to obtain
non-vanishing terms
$$\Fbar_{k+1}(x,y)
\;\;=\;\;\sum_{m<r}\;\;
\Fbar_k
(s_{[m,m+n-r]}\,x_{0\ldots m},\,
y_{0 \ldots m r \ldots n}
)
\;\otimes\;
x_{m \ldots r}$$
The theorem then follows by a straightforward induction.
\qed\enddemo

\demo{Proof of Lemma 5.3}
Expanding the definitions (5.1--5.3) we get
$$\align
\varphi_n(x,y)&=\qquad \sum_{m=0}^{n-1} \quad\pm g_{n-m}^{(m+1)}s_m\;(x,y)\\
&=
\sum\Sb 0\le m\le n-1\\ 0\le \ell \le n-m\\ A\cup B=[0,n-m-1] \\
    |A|=n-m-\ell ,|B|=\ell \endSb \pm
(s_{A+m+1}\del _{\ell +1+m+1}^{\;n-m-\ell }s_m\,x
,s_{B+m+1}\del _{0+m+1}^{\;\ell } s_m\,y)\\
&=\sum\Sb0\le m\le n-1\\ 1\le \ell \le n-m\\ A\cup B=[m+1,n] \\
    |A|=n-m-\ell ,|B|=\ell \endSb \pm
(s_As_m\del _{m+\ell +1}^{\;n-m-\ell }\,x
,s_B\del _{m+1}^{\;\ell -1}\,y)\\
&=\sum\Sb 0\le m\le n-1\\ m+1\le r\le n\\ A\cup B=[m+1,n] \\
    |A|=n-r,|B|=r-m\endSb \pm
(s_{A\cup\{m\}}\,x_{0\ldots r},\,s_B\,y_{0\ldots m r\ldots n}).
\endalign$$
Here we have dropped the $\ell =0$ terms, since they are degenerate (in the
image of $s_m$), and written $r=m+\ell $.
\qed\enddemo

\demo{Proof of Lemma 5.4}
By Lemma 5.2 the first equation holds for $\ell\le m$, since then
$$\Fbar_1((\varphi^{A,B}(x,y))_{0\ldots\ell}
)\;\;=\;\;\Fbar_1((x,y)_{0\ldots\ell})
\;\;=\;\;x_{0\ldots\ell}\;+\;y_{0\ldots\ell}$$
If $\ell>m$, then the first term is always $s_m$-degenerate,
as is the second, unless all the indices specified in $B$ are
$\geq\ell$, that is, unless $[m+1,\ell-1]\subseteq A$. 

If $\ell\le m$, then $(\varphi^{A,B}(x,y))_{\ell\ldots n+1}$
always has  an $s_m$-degeneracy in the first component and some
other degeneracy in the second, since $B\neq\emptyset$, so its image
under $\Fbar_1$ is zero. If $\ell>m$, then the first component can
only be non-degenerate if no elements of $A$ are $\geq\ell$, so
$[\ell,n]\subseteq B$. 
Similarly the second can be non-degenerate only if $[\ell,n]\subseteq A$. 

If $\ell>m$ then every term of $\varphi((\varphi^{A,B}(x,y))_{0,\ldots,\ell})$
has as first factor the face of the second component specified by
$[0,p]\cup[q,\ell]$, say, which
is degenerate if $q\leq m$ and $[q,\ell-1]\cap B$ is
non-empty. If $q>m$, however, then
either $p>m$ and $[p,q]\subseteq B$, so it is still degenerate, or
$[p,q]\cap(A\cup\{m\})$ is non-empty and the other factor is
degenerate.

A similar argument shows that the terms of
$\varphi((\varphi^{A,B}(x,y))_{\ell\ldots n+1})$ 
always have one of the two factors degenerate.
\qed\enddemo


\enddocument